\input amstex
\documentstyle{amsppt}
\nologo
\magnification=1200
\pageheight{19cm}

\loadbold

\define\inn#1{{\overset{\lower1.5pt\hbox{$\ssize\circ$}}\to#1}}

\define\eps{\varepsilon}
\redefine\phi{\varphi}

\define\w{{\frak w}}
\define\LL{{\Cal L}}

\define\RR{{\Cal R}}
\define\J{{\Cal J}}

\define\h{{\frak z}}

\define\HH{{\Cal H}}

\define\la{\lambda}
\define\lap{\lambda'}

\define\subla{_\lambda}

\define\sublap{_{\lambda'}}

\define\om{\omega}
\define\Om{\Omega}
\define\gfs{{g_{\text{FS}}}}

\define\scp{{\<\,.\,,.\,\>}}
\define\inv{^{-1}}
\define\Sym{\text{\rm{Sym}}}

\define\Aut{{\text{\rm{Aut}}}}
\define\Autm{\Aut{}_{g_0}^{\,T}(M)}
\define\Autbm{\overline{\Aut}{}_{g_0}^{\,\,T}(M)}

\define\Id{{\text{\rm{Id}}}}
\define\SU{{\text{\rm{SU}}}}
\define\su{{\text{$\frak s$}\text{$\frak u$}}}

\redefine\O{{\text{\rm{O}}}}
\define\SO{{\text{\rm{SO}}}}
\define\so{{\text{$\frak s$}\text{$\frak o$}}}

\define\dvol{{\text{{\it dvol}}}}

\define\spann{{\text{\rm{span}}\,}}
\define\kernn{{\text{\rm{ker}}\,}}

\define\dimm{{\text{\rm{dim}}}}

\redefine\exp{{\text{\rm{exp}}}}

\define\tr{{\text{\rm{tr}}}}
\define\scal{{\text{\rm{scal}}}}

\define\Isom{{\text{\rm{Isom}}}}

\define\N{{\Bbb N}}
\define\R{{\Bbb R}}
\define\C{{\Bbb C}}
\define\Z{{\Bbb Z}}
\define\<{{\langle}}
\define\>{{\rangle}}
\define\minzero{\setminus\{0\}}

\define\restr#1{{\lower0.4ex\hbox{$\vert$}\lower0.7ex
  \hbox{$\ssize{#1}$}}}
\define\drestr#1{{\lower0.6ex\hbox{$\vert$}\lower1ex
  \hbox{$\ssize{#1}$}}}
\define\tdiffzero#1{{\frac d{d{#1}}\lower0.2ex\hbox{\restr{{#1}=0}}\,}}
\define\diffzerosec#1{{\frac{d^2}{d{#1}^2}\lower1ex\hbox{$\big
  \vert_{\raise1\jot\hbox{${\tsize {#1}=0}$}}$}\,}}

\topmatter

\title Isospectral metrics on five-dimensional spheres
\endtitle

\rightheadtext{Isospectral metrics on five-dimensional spheres}

\author Dorothee Schueth\endauthor

\address Mathematisches Institut, Universit\"at Bonn, Beringstr\. 1,
         D@-53115 Bonn, Germany\endaddress

\email schueth\@math.uni-bonn.de\endemail

\thanks The author is partially supported by SFB~256, Bonn.
\endthanks


\abstract
We construct isospectral pairs of Riemannian metrics on $S^5$
and on $B^6$, thus lowering by three the minimal dimension of spheres and balls
on which such metrics have been constructed previously ($S^{n\ge8}$ and
$B^{n\ge9}$). We also construct
continuous families of isospectral Riemannian metrics on~$S^7$ and on~$B^8$.
In each of these examples, the metrics can be chosen equal to the
standard metric outside subsets of arbitrarily small volume.
\endabstract

\endtopmatter

\document

\heading Introduction\endheading

\noindent
During the past two decades, research on isospectral manifolds --
that is, Riemannian manifolds sharing the same spectrum (including
multiplicities) of the Laplace operator acting on functions --
has been very active; see, for example, the survey article~\cite{Go2}.
However, it was only recently that Zoltan I.~Szab\'o
and Carolyn Gordon independently discovered the first examples
of isospectral metrics on spheres: Pairs of such metrics on~$S^{n\ge10}$
\cite{Sz2,3}, and continuous families on~$S^{n\ge8}$ \cite{Go3}.

In spite of the wealth of other isospectral manifolds obtained before,
the construction of isospectral spheres had seemed beyond reach
for a long time. The~reason was that both of the main
methods of construction which were known and used until~2000 had
excluded spheres:
\roster
\item"$\bullet$"
The so-called Sunada method~\cite{Su} and its various
generalizations (see, e.g.,~\cite{DG}) produces isospectral quotients
$(M/\varGamma_1\,,g)$, $(M/\varGamma_2\,,g)$ of a common Riemannian covering
manifold $(M,g)$; in particular, these isospectral manifolds were
always nonsimply connected. Many interesting examples of locally
isometric isospectral manifolds of the form $(M/\varGamma_1\,,g)$,
$(M/\varGamma_2\,,g)$ -- some of Sunada type, some using other special
constructions -- can be found, for example, in~\cite{Vi}, \cite{Ik},
\cite{GW1}, \cite{DG}, \cite{Gt1,2}, \cite{Sch1}; also the famous
examples of isospectral plane domains~\cite{GWW} arise as orbifolds
in the Sunada type setting.
\item"$\bullet$"
The method of principal torus bundles (see~\cite{Go1}, \cite{GW2},
\cite{GGSWW}, \cite{Sz1}, \cite{GSz}, \cite{Sch2,3}) produces certain
pairs of isospectral principal bundles for which the structural group
is a torus of dimension at least two; the metrics are invariant
under the torus action. Since no sphere is a principal
$T^{k\ge2}$@-bundle, spheres cannot be obtained by using this method
either (although products of spheres could~\cite{Sch2}; these
were the first examples of simply connected isospectral manifolds).
See~\cite{Sch3} for a detailed treatment of the method
of principal torus bundles, a systematical approach for applying it,
and many examples (among others, isospectral left invariant metrics
on compact Lie groups).
\endroster

The key to Gordon's construction of isospectral metrics on
the spheres $S^{n\ge8}$ (and on the balls~$B^{n\ge9}$) was a new
approach which still involves $T^{k\ge2}$@-actions, but does not require
them to be free anymore. On the other hand, Z.I.~Szab\'o's
pairs of isospectral metrics on~$S^{n\ge10}$ (and on~$B^{n\ge11}$) do not
arise by such a construction and seem to be of a completely different
type. A spectacular feature of his examples is that they include
pairs of isospectral metrics on~$S^{11}$ in which one of the metrics
is homogeneous while the other is not.

\smallskip

The present paper serves several purposes.

First, we reformulate Gordon's new theorem \cite{Go3, Theorem~1.2}
in a somewhat more elegant way -- see Theorem~1.4 below --, and
in Theorem~1.6 we establish a special version of it which turns
to be a useful tool for finding new applications. It~also accounts
for all applications of Theorem~1.4 which are known so far.

Second, we derive a general sufficient nonisometry condition for
the type of Riemannian metrics occurring in Theorem~1.6;
see Proposition~2.4
below. Although we will apply this criterion only to certain new examples
constructed in this paper, we wish to point out that it could also be used
to unify most of the various nonisometry proofs from \cite{GW2},
\cite{GGSWW}, \cite{GSz}, \cite{Go3} (however, some of the locally
homogeneous manifolds with boundary from \cite{GW2}
and~\cite{GSz} cannot be proven nonisometric using  this approach
because they violate
a certain genericity condition; see~\thetag{G} in Proposition~2.3).
Similarly, it would be possible
to apply our method to the isospectral examples from \cite{Sch2,3} where
we had instead performed explicit curvature computations in
order to show that objects like $\int\scal^2$, $\int\|R\|^2$,
the critical values of the scalar curvature, or the dimension
of their loci are not spectrally determined.

Our third (and main) purpose is to use Theorem~1.6 for constructing
isospectral pairs of metrics on~$S^5$ (and on~$B^6$) and thus decreasing
the minimal dimensions of Gordon's examples by three; see Example~3.3.
We also obtain continuous isospectral families of metrics on~$S^7$
(and on~$B^8$); see Example~3.2.
For the nonisometry proofs we use the general criterion mentioned above.

Fourth, we show that in these new examples -- in particular, on~$S^5$
and~$B^6$
(pairs), and on~$S^7$ and~$B^8$ (continuous families)
-- it is possible to choose the
isospectral metrics in such a way that they are equal to the round
(resp\. flat) metric outside certain subsets of arbitrarily small volume;
see Theorm~5.3.
This gives a nice contrast to the fact that in dimensions up to six,
the round spheres themselves
are completely determined by their spectra~\cite{Ta1} (the corresponding
problem in higher dimensions has not yet been solved), and to the fact
that no round sphere in any dimension admits a continuous isospectral
deformation~\cite{Ta2}.

\smallskip

The paper is organized as follows:

In Section~1 we present our reformulation (Theorem~1.4) of Gordon's
theorem and the afore-mentioned specialization (Theorem~1.6).
Section~2 contains our general sufficient nonisometry criterion
for the metrics occurring in Theorem~1.6.
In~Section~3 we construct our new examples: Continuous families
of isospectral metrics on~$S^7$ and on~$B^8$ (Example~3.2), and
pairs of such metrics on~$S^5$ and on~$B^6$ (Examle~3.3).
Moreover, we give a survey of a number of related examples in~3.4;
each of them can be obtained using Theorem~1.6.
Section 4 gives the nonisometry proof for Examples~3.2\,/\,3.3.
Finally, we show in Section~5
how to make the isospectral metrics round (resp\. flat)
on large subsets (Theorem~5.3).

\smallskip

The author would like to thank Carolyn Gordon and Werner Ballmann
for interesting discussions concerning these and related topics.

\bigskip

\heading\S1 Isospectrality via effective torus actions\endheading

\definition{1.1 Definition}
(i) The {\it spectrum\/} of a closed Riemannian manifold is the spectrum
of eigenvalues, counted with multiplicities, of the associated Laplace
operator acting on functions. The {\it Dirichlet spectrum\/} of a compact
Riemannian manifold~$M$ with boundary is the spectrum of eigenvalues
corresponding to eigenfunctions which satisfy the Dirichlet boundary
condition $f\restr{\partial M}=0$.
The {\it Neumann spectrum} of such a manifold
is defined analogously with respect to the Neumann boundary condition
$Nf=0$, where~$N$ is the inward-pointing unit normal field on the boundary.

(ii) Two closed Riemannian manifolds are called {\it isospectral\/}
if they have the same spectrum (including multiplicities). Two compact
Riemannian manifolds with boundary are called {\it Dirichlet isospectral},
resp\. {\it Neumann isospectral}, if they have the same Dirichlet spectrum,
resp\. the same Neumann spectrum.
\enddefinition

\subheading{1.2 Remark}
Let $(M,g)$ be a compact Riemannian manifold with or without boun\-dary.
Consider the Hilbert spaces
$\HH^N:=H^{1,2}(M,g)$ and $\HH^D:=\inn{H}{}^{1,2}(M,g)\subseteq\HH^N$.
Note that~$\HH^N$, resp\. $\HH^D$, is the completion of
$C^\infty(M)$, resp\. $\{f\in C^\infty(M)\mid f\restr{\partial M}=0\}$,
with respect to the $H^{1,2}$@-norm associated to~$(M,g)$.
For each $f\in\HH^N\minzero$ the Rayleigh quotient is defined as
$$\RR(f):=\tsize\int_M\|df\|_g^2\dvol_g\Big/\tsize\int_M|f|^2\dvol_g=
\bigl(\|f\|^2_{H^{1,2}(M,g)}\bigm/\|f\|^2_{L^2(M,g)}\bigr)-1.
$$
Let $0\le\la_1^N\le\la_2^N\le\ldots\to\infty$, resp\.
$0\le\la_1^D\le\la_2^D\le\ldots\to\infty$, denote the Neumann spectrum,
resp\. the Dirichlet spectrum, of $(M,g)$; if $\partial M=\emptyset$
then both of these sequences coincide with the spectrum of $(M,g)$.
Finally, denote by $L_k^N$, resp.~$L_k^D$, the set of all $k$@-dimensional
subspaces of~$\HH^N$, resp.~$\HH^D$. Then we have the following
{\it variational characterization of eigenvalues\/} (see, e.g.,~\cite{Be}):
$$\la_k^N=\inf_{U\in L_k^N}\sup_{f\in U\minzero}\RR(f)\quad\text{ and }\quad
  \la_k^D=\inf_{U\in L_k^D}\sup_{f\in U\minzero}\RR(f).\tag{1}  
$$  

\subheading{1.3 Notation}
By a {\it torus\/}, we always mean a nontrivial, compact, connected abelian
Lie group. If a torus~$T$ acts smoothly and effectively by isometries
on a compact connected
Riemannian manifold~$(M,g)$ then we denote by~$\hat M$ the union of those
orbits on which~$T$ acts freely.
Note that~$\hat M$ is an open dense submanifold of~$M$.
The~action of~$T$ gives~$\hat M$ the
structure of a principal $T$@-bundle. By~$g^T$ we denote the unique
Riemannian
metric on the quotient manifold $\hat M/T$ such that the canonical
projection $\pi:(\hat M,g)\to (\hat M/T,g^T)$ is a Riemannian
submersion.

\proclaim{1.4 Theorem}
Let $T$ be a torus which acts effectively on two compact connected
Riemannian manifolds $(M,g)$ and $(M',g')$ by isometries.
For each subtorus~$W\subset T$ of codimension one, suppose that there exists
a $T$@-equivariant diffeomorphism $F_W: M\to M'$ which
satisfies $F_W^*\dvol_{g'}=\dvol_{g}$
and induces an isometry $\bar F_W$ between the quotient manifolds
$(\hat M/W,g^W)$ and $(\hat M'/W,g^{\prime\,W})$.
Then $(M,g)$ and $(M',g')$ are isospectral; if the manifolds have
boundary then they are Dirichlet and Neumann isospectral.
\endproclaim

\remark{Remark}
Theorem 1.4 above is a slight variation of Carolyn Gordon's Theorem~1.2
in~\cite{Go3}. Instead of our condition $F_W^*\dvol_{g'}=\dvol_{g}$\,,
Gordon assumes the condition that $\bar F_{W*}$ maps the projected
mean curvature vector field~$\bar H_W$ of the submersion
$(\hat M,g)\to(\hat M/W,g^W)$ to the corresponding vector
field~$\bar H'_W$\,. While the two conditions actually turn out to be
equivalent in this context, our volume preserving condition is not only
easier to formulate but also more convenient to check in applications.
More than that, it is inherent and automatically satisfied in the
specialization described below in~1.5\,/\,1.6 which covers all
applications of the above theorem which are known so far.
Our different formulation of the theorem has also led to a different
proof (construction of a certain isometry between the $H^{1,2}$@-spaces
instead of intertwining the Laplacians). An advantage of our proof is that
it does not require a certain additional condition which Gordon assumes
in the case of manifolds with boundary (namely, that $\hat M\cap \partial M$
be dense in~$\partial M$).
\endremark

\demo{Proof of Theorem 1.4}
Consider the Hilbert space $\HH:=H^{1,2}(M,g)$
in the case of manifolds without boundary or in case of Neumann
boundary conditions, resp\. $\HH:=\inn{H}{}^{1,2}(M,g)$
in the case of Dirichlet boundary conditions.
Let $\HH'$ be defined analogously with respect to $(M',g')$.
We~claim that there is a Hilbert space isometry from $\HH$ to $\HH'$
which moreover preserves $L^2$@-norms;
the theorem will thus follow from the variational characterization of
eigenvalues~\thetag{1}.

Consider the unitary representation of~$T$ on~$\HH$ defined by
$(zf)(x)=f(zx)$ for all $f\in\HH$, $z\in T$, $x\in M$.
Write $T=\h/\LL$ and let $\LL^*$ be the dual lattice.
Since~$T$ is abelian, $\HH$ decomposes as the orthogonal sum
$\bigoplus_{\mu\in\LL^*}\HH_\mu$ with
$\HH_\mu=\{f\in\HH\mid zf=e^{2\pi i\mu(Z)}f\text{ \,for
all }z\in T\}$,
where~$Z$ denotes any representative for~$z$ in~$\h$.
In particular, this implies the coarser decomposition
$$\HH=\HH_0\oplus{\tsize\bigoplus}_W(\HH_W\ominus\HH_0),\tag{2}
$$
where $W$ runs though the set of all subtori of codimension~$1$ in~$T$,
and $\HH_W$ is the sum of all $\HH_\mu$ such that $\mu\in\LL^*$
and $T_eW\subseteq\kernn\mu$. In other words, $\HH_W$ is just the space
of $W$@-invariant functions in~$\HH$.
Let $\HH'_W$ and~$\HH'_0$ be the analogously defined subspaces of~$\HH'$.
Now let~$W$ be any subtorus of codimension~$1$ in~$T$ and choose a
diffeomorphism~$F_W$ as in the assumption. Since~$F_W$ intertwines the
$T$@-actions, $F_W^*$ sends~$\HH'_W$
to~$\HH_W$ and $\HH'_0\subset\HH'_W$ to $\HH_0\subset
\HH_W$\,. We will now show that $F_W^*:\HH'_W\to\HH_W$ is a Hilbert space
isometry which also preserves $L^2$@-norms; in view of the
decomposition~\thetag{2} this will prove our
above claim. Preservation of $L^2$@-norms is trivial by the assumption
$F_W^*\dvol_{g'}=\dvol_g$\,. Moreover, this assumption implies that
it suffices to show that for each $\psi\in C^\infty(M')$ which is
invariant under  the $W$@-action and for all $y\in M'$ we have
$\|d\psi\restr y\|_{g'}=\|d\phi\restr x\|_g$\,,
where $\phi:=F_W^*\psi$ and $x:=F_W\inv(y)$.
We can assume $x\in\hat M$; let~$\bar\phi$ and~$\bar\psi$
be the functions induced
on $\hat M/W$ and $\hat M'/W$. Then $\bar\phi=\bar F_W^*\bar\psi$\,;
since $g^W$ and~$g^{\prime\,W}$
are the submersion metrics and~$\bar F_W$ is an isometry
we obtain indeed, at the appropriate points: $\|d\phi\|_g
=\|d\bar\phi\|_{g^W} = \|d\bar\psi\|_{g^{\prime\,W}} = \|d\psi\|_{g'}$\,.
\qed\enddemo

\subheading{1.5 Notation and Remarks}
In the following we fix a torus~$T$ with Lie algebra~$\h=T_eT$.
Let~$\LL$ be the cocompact lattice in~$\h$ such that $\exp:\h\to T$
induces an isomorphism from $\h/\LL$ to~$T$, and denote by~$\LL^*
\subset\h^*$ the dual lattice. We also fix a compact connected Riemannian
manifold
$(M,g_0)$, with or without boundary, and a smooth effective action of~$T$
on $(M,g_0)$ by isometries.
\roster
\item"(i)"
For $Z\in\h$ we denote by $Z^*$ the vector field
$x\mapsto\tdiffzero t\,\exp(tZ)x$ on~$M$.
For each $x\in M$ and each subspace~$\w$ of~$\h$ we let
$\w_x:=\{Z^*_x\mid Z\in\w\}$.
\item"(ii)"
We call a smooth $\h$@-valued $1$@-form on~$M$ {\it admissible\/}
if it is $T$@-invariant and horizontal (i.e., vanishes
on the vertical spaces~$\h_x$).
\item"(iii)"
For any admissible $\h$@-valued $1$@-form~$\la$ on~$M$ we denote
by~$g\subla$ the Riemannian metric on~$M$ given by
$$g\subla(X,Y)=g_0\bigl(X+\la(X)^*,Y+\la(Y)^*\bigr).
$$
In other words,
$g\subla=(\Phi\subla\inv)^*g_0$\,, where~$\Phi\subla$ is the smooth
endomorphism field on~$M$ given by $X\mapsto X-\la(X)^*$ for all $X\in TM$.
Note that~$\Phi_\la$ is unipotent on each tangent space;
in particular, $\dvol_{g\subla}=\dvol_{g_0}$\,.
\item"(iv)"
Finally, note that~$g\subla$ is again invariant under the action
of~$T$, that~$g\subla$ restricts to the same metric as~$g_0$ on
the vertical subspaces~$\h_x$\,, and that the submersion metric~$g\subla^T$
on~$\hat M/T$ is equal to~$g_0^T$\,.
\endroster

\proclaim{1.6 Theorem}
In the context of Notation~{\rm 1.5},
let $\la$, $\lap$ be two admissible $\h$@-valued $1$@-forms on~$M$.
Assume:
\roster
\item"($*$)" For every $\mu\in\LL^*$ there exists a $T$@-equivariant
  $F_\mu\in\Isom(M,g_0)$ which satisfies $\mu\circ\la=F_\mu^*(\mu\circ\lap)$.
\endroster
Then $(M,g\subla)$ and $(M,g\sublap)$ are isospectral; if~$M$ has boundary
then the two manifolds are Dirichlet and Neumann isospectral.
\endproclaim

\demo{Proof}
We show that $(M,g\subla)$ and $(M,g\sublap)$ satisfy the hypotheses
of Theorem~1.4. Let~$W$ be a subtorus of codimension~$1$ in~$T$.
Choose $\mu\in\LL^*$ such that $\w:=T_eW$ equals $\kernn\mu$,
and choose a corresponding~$F_\mu$ as in~\thetag{$*$}.
We claim that $F_W:=F_\mu$ satisfies the conditions required in Theorem~1.4.
First of all note that~$F_\mu$\,, being an isometry of~$g_0$\,, trivially
satisfies $F_\mu^*\dvol_{g\sublap}=\dvol_{g\subla}$ because of
$\dvol_{g\subla}=\dvol_{g_0}=\dvol_{g\sublap}$ (see 1.5(iii)).
Thus it remains to prove that~$F_\mu$ induces an isometry from
$(\hat M/W,g\subla^W)$ to $(\hat M/W,g\sublap^W)$.

Let~$V\in T_xM$ be any vector which is $g\subla$@-orthogonal to~$\w_x$\,;
then $V=\Phi\subla(X)$ for some~$X$ which is
$g_0$@-orthogonal to~$\w_x$\,. By condition~\thetag{$*$} we know that
$\la(X)$ equals $\lap(F_{\mu*}X)$ modulo~$\w$. Keeping in mind
that~$F_\mu$ commutes with the~$T$@-action, we conclude that the vector
$F_{\mu*}V=F_{\mu*}(\Phi\subla X)$ equals $Y:=\Phi\sublap(F_{\mu*}X)$ up to
an error in $\w_{F_\mu(x)}$\,. But $F_{\mu*}X$ is $g_0$@-orthogonal to
$\w_{F_\mu(x)}$\,; thus $Y$ is the projection of $F_{\mu*}V$ to the
$g\sublap$@-orthogonal complement of~$\w_{F_\mu(x)}$\,.
Our assertion now follows from $\|Y\|_{g\sublap}=\|F_{\mu*}X\|_{g_0}
=\|X\|_{g_0}=\|V\|_{g\subla}$\,.
\qed\enddemo   

\bigskip

\heading\S2 A sufficient condition for nonisometry\endheading

\noindent
Throughout this section we let $(M,g_0)$, $T$, $\h$ be as in Notation~1.5,
and we define the principal $T$@-bundle~$\pi:\hat M\to\hat M/T$
as in Notation~1.3.
By $\la,\lap$ we will always denote admissible
$\h$@-valued $1$@-forms on~$M$.
In Proposition~2.4 below we will establish a sufficient condition
for $(M,g\subla)$, $(M,g\sublap)$ to be nonisometric.

\subheading{2.1 Notation and Remarks}
\roster
\item"(i)"
We say that a diffeomorphism $F:M\to M$ is {\it $T$@-preserving\/}
if conjugation by~$F$ preserves $T\subset\text{Diffeo(M)}$.
In that case, we denote by~$\Psi_F$ the automorphism of~$\h=T_eT$
induced by conjugation by~$F$. Obviously, each $T$@-preserving
diffeomorphism~$F$ of~$M$ maps $T$@-orbits to $T$@-orbits and satisfies
$F_*(Z^*)=\Psi_F(Z)^*$ for all $Z\in\h$, where the vector fields~$Z^*$
on~$M$ are defined as in~1.5(iii).
\item"(ii)"
We denote by $\Autm$ the group of all $T$@-preserving diffeomorphisms~$F$
of~$M$ which, in addition, preserve the $g_0$@-norm of vectors tangent
to the $T$@-orbits and induce an isometry of $(\hat M/T,g_0^T)$.
We denote the corresponding group of induced isometries by
$\Autbm\subset\Isom(\hat M/T,g_0^T)$.
\item"(iii)"
We define $\Cal D:=\{\Psi_F\mid F\in\Autm\}
\subset\Aut(\h)$.
Note that~$\Cal D$ is discrete because it is a subgroup of the discrete
group $\{\Psi\in\Aut(\h)\mid\Psi(\LL)=\LL\}$, where~$\LL$ is the lattice
$\kernn(\exp:\h\to T)$.
\item"(iv)"
Let~$\om_0$ denote the connection form on~$\hat M$ associated
with~$g_0$\,; i.e., for each $x\in\hat M$ the horizontal space
$\kernn(\om_0\restr{T_x\hat M})$ is the $g_0$@-orthogonal
complement of~$\h_x$ in~$T_x\hat M$. Then
the connection form on~$\hat M$ associated with~$g\subla$ 
is obviously given by $\om\subla:=\om_0+\la$.
\item"(v)"
Let~$\Om\subla$ denote the curvature form on $\hat M/T$ associated
with the connection form~$\om\subla$ on~$\hat M$.
We have $\pi^*\Om\subla=d\om\subla$ because $T$~is abelian.
\item"(vi)"
Since~$\la$ is $T$@-invariant and horizontal it induces some~$\h$@-valued
$1$@-form~$\bar\la$ on~$\hat M/T$. We conclude from $\pi^*\Om\subla=
d\om\subla=d\om_0+d\la$ that $\Om\subla=\Om_0+d\bar\la$. In particular,
$\Om\subla$~and~$\Om_0$ differ by an exact $\h$@-valued $2$@-form.
\endroster

\proclaim{2.2 Lemma}
Suppose that $F:(M,g\subla)\to(M,g\sublap)$ is a $T$@-preserving isometry.
\roster
\item"(i)"
$F$ preserves the $g_0$@-norm of vectors tangent to the $T$@-orbits,
and it induces an isometry~$\bar F$ of $(\hat M/T,g_0^T)$.
In particular, $F\in\Autm$ and $\Psi_F\in\Cal D$.
\item"(ii)"
$F^*\om\sublap=\Psi_F\circ\om\subla$\,; in particular,
$F^*d\om\sublap=\Psi_F\circ d\om\subla$\,.
\item"(iii)"
The isometry~$\bar F$ of $(\hat M/T,g_0^T)$ satisfies 
$$\bar F^*\Omega\sublap=\Psi_F\circ\Omega\subla\,.\tag{3}
$$
\endroster
\endproclaim

\demo{Proof}
(i) This follows immediately from~1.5(iv).

(ii) We have $\om\sublap(Z^*)=\om\subla(Z^*)=Z$ for all $Z\in\h$,
hence $\om\sublap(F_*(Z^*))=\om\sublap(\Psi_F(Z)^*)=\Psi_F(Z)=
\Psi_F(\om\subla(Z^*))$. The equation thus holds when applied to
vectors tangent to the $T$@-orbits.
Since~$F$ is an isometry and maps orbits to orbits, it
must map $g\subla$@-horizontal vectors to
$g\sublap$@-horizontal
vectors. Hence both sides of the asserted equation vanish when applied
to a $g\subla$@-horizontal vector.

(iii)
This follows from~(ii) and~2.1(v).
\qed\enddemo

\proclaim{2.3 Proposition}
Let~$\la$ be an admissible $\h$@-valued $1$@-form on~$M$ such that the
associated curvature form~$\Om\subla$ on~$\hat M/T$ satisfies the
following genericity condition:
$$\text{No nontrivial ${}\,1$@-parameter group in ${}\,\Autbm$
preserves~$\Om_\la$\,.}
\tag{G}
$$
Then~$T$ is a maximal torus in $\Isom(M,g\subla)$.
\endproclaim

\demo{Proof}
Let $F_t\in\Isom(M,g\subla)$ be a $1$@-parameter family of isometries
commuting with~$T$. In particular, the maps~$F_t$ are $T$@-preserving
and thus induce a $1$@-parameter family $\bar F_t
\in\Isom(\hat M/T,g_0^T)$.
The corresponding $\Psi_{F_t}\in\Cal D$ satisfy $\Psi_{F_t}\equiv\Id$
because~$\Cal D$ is discrete and $\Psi_{F_0}=\Id$.
Thus according to~\thetag{3}, each~$\bar F_t$ preserves~$\Om\subla$\,.
The assumed property~\thetag{G} of~$\Om\subla$ now implies
$\bar F_t\equiv\Id$. We conclude that~$F_t$
restricts to a gauge transformation of the principal $T$@-bundle~$\hat M$.
On the other hand, $F_t^*\om\subla\equiv\om\subla$ by~2.2(ii).
But a gauge transformation which preserves the connection form
of a principal bundle must act as an element of the
structural group on each connected component of the bundle.
Since an isometry of the connected Riemannian manifold $(M,g\subla)$
is determined by its restriction to any nonempty open subset, it
follows that the family~$F_t$ is contained in~$T$. 
\qed\enddemo

\proclaim{2.4 Proposition}
Let $\la,\lap$ be admissible $1$@-forms on~$M$ such that~$\Om\sublap$
has property~\thetag{G}.
Furthermore, assume that
$$\Om\subla\notin\Cal D\circ\Autbm^*\Om\sublap\,.\tag{N}
$$
Then $(M,g\subla)$ and $(M,g\sublap)$ are not isometric.
\endproclaim

\demo{Proof}
Suppose that there were an isometry $F:(M,g\subla)\to(M,g\sublap)$.
By Proposition~2.3, $T$~is a maximal torus in $\Isom(M,g\sublap)$.
Since all maximal tori are conjugate, we can assume~$F$ -- after possibly
combining it with an isometry of~$(M,g\sublap)$ -- to be $T$@-preserving.
But then Lemma~2.2 implies $\bar F^*\Om\sublap=\Psi_F\circ\Om\subla$
with $\bar F\in\Autbm$ and $\Psi_F\in\Cal D$,
which contradicts our assumption~\thetag{N}.
\qed\enddemo

\subheading{2.5 Remark}
Note that equally valid (but weaker)
versions of Propositions~2.3 and~2.4
would be obtained by replacing $\Autbm$ by the possibly larger group
$\Isom(\hat M/T,g_0^T)$.

\bigskip

\heading\S3 Examples\endheading

\subheading{3.1 Notation}
\noindent
\roster
\item"(i)"
Throughout the new examples given in~3.2 and~3.3 below
we consider the two-dimensional
torus $T:=\R^2/\LL$ with $\LL:=2\pi\Z\times 2\pi\Z$.
We denote the standard basis of its Lie algebra~$\h\cong\R^2$ by
$\{Z_1\,,Z_2\}$.
\item"(ii)"
We let $\Cal E$ be the group of the four linear isomorphisms
of~$\h$ which preserve each of the sets $\{\pm Z_1\}$ and $\{\pm Z_2\}$.
\item"(iii)"
Let $m\in\N$. We identify the real vector spaces $\C^{m+1}=\C^m\oplus\C$
and~$\R^{2m+2}$ via the linear isomorphism which sends
$\{e_1\,,ie_1\,,\,\ldots,e_{m+1}\,,ie_{m+1}\}$ (in this order) to
the standard basis of~$\R^{2m+2}$, where $\{e_1\,,\,\ldots,e_{m+1}\}$
denotes the standard basis of~$\C^{m+1}$.
We let the torus~$T$ act on this space by
$$\exp(aZ_1+bZ_2):(p,q)\mapsto(e^{ia}p,e^{ib}q)
$$
for all $a,b\in\R$, $p\in\C^m$, $q\in\C$.
This action preserves the unit sphere $S^{2m+1}\subset\R^{2m+2}$ as well
as the unit ball $B^{2m+2}$;
by restriction we thus obtain an action~$\rho$ of~$T$ on~$S^{2m+1}$,
respectively on~$B^{2m+2}$.
\endroster

\subheading{3.2 Example:\newline
Continuous isospectral families of metrics on~$\boldkey S^{\boldkey 2\boldkey m
\boldkey+\boldkey 1\boldsymbol\ge\boldkey 7}$
and on $\boldkey B^{\boldkey 2\boldkey m\boldkey+
\boldkey 2\boldsymbol\ge\boldkey 8}$}

\noindent
\subsubhead 3.2.1 Notation\endsubsubhead
For each linear map $j:\h\cong\R^2\to\su(m)$ we define a $\h$@-valued
$1$@-form $\la=(\la^1,\la^2)$ on $\R^{2m+2}\cong\C^m\oplus\C$ by letting
$$\la^k_{(p,q)}(X,U)=|p|^2\<j_{Z_k}p,X\>-\<X,ip\>\<j_{Z_k}p,ip\>\tag{4}
$$
for $k=1,2$ and all $(X,U)\in T_p\R^{2m}\oplus T_q\R^2$.
Here $\scp$ denotes the standard euclidean inner product on~$\R^{2m}$,
and the skew-hermitian maps $j_{Z_k}:=j(Z_k)$ act on $p\in\R^{2m}$
via the above identification $\R^{2m}\cong\C^m$.

By restriction we obtain a smooth $\h$@-valued $1$@-form~$\la$ on the
unit sphere~$S^{2m+1}$, respectively on the unit ball~$B^{2m+2}$.

\subsubhead 3.2.2 Remark\endsubsubhead
We observe that~$\la$ is admissible (see Notation~1.5(ii))
with respect to the action~$\rho$ of~$T$ defined above in~3.1(iii).
In fact, invariance of~$\la$ under the action of~$T$ is immediate because
multiplication with the complex scalar factor~$e^{ia}$ commutes
with~$j_{Z_k}\in\su(m)$ and preserves the euclidean inner product.
It remains to check that~$\la$ vanishes on the spaces~$\h_{(p,q)}
=\spann\{(ip,0),(0,iq)\}\subset T_p\R^{2m}\oplus T_q\R^2$. Indeed we have
$\la^k_{(p,q)}(ip,0)=|p|^2\<j_{Z_k}p,ip\>-\<ip,ip\>\<j_{Z_k}p,ip\>=0$
for $k=1,2$, and $\la_{(p,q)}(0,iq)=0$ by definition.

\subsubhead 3.2.3 Definition\endsubsubhead
Let~$g_0$ be the round standard metric on~$S^{2m+1}$, respectively
the standard metric on~$B^{2m+2}$,
and let~$g\subla$ be the metric associated with~$\la$ and~$g_0$ as
in Notation~1.5(iii).

\smallskip

Summarizing, for each linear map $j:\R^2\to\su(m)$ we have an associated
Riemannian metric~$g\subla$ on~$S^{2m+1}$, respectively on $B^{2m+2}$,
via the corresponding $1$@-form~$\la$ as defined in~\thetag{4}.

\subsubhead 3.2.4 Definition\endsubsubhead
Let $j,j':\h\cong\R^2\to\su(m)$ be two linear maps.
\roster
\item"(i)"
We call $j$ and~$j'$ {\it isospectral}, denoted $j\sim j'$,
if for each $Z\in\h$ there exists $A_Z\in\SU(m)$ such that
$j'_Z=A_Zj_ZA_Z\inv$\,.
\item"(ii)"
Let $Q:\C^m\to\C^m$ denote complex conjugation.
We call $j$ and~$j'$ {\it equivalent}, denoted $j\cong j'$,
if there exists $A\in\SU(m)\cup\SU(m)\circ Q$
and $\Psi\in\Cal E$ (see Notation~3.1(ii))
such that $j'_Z=Aj_{\Psi(Z)}A\inv$ for all $Z\in\h$.
\item"(iii)" We say $j$ is {\it generic\/} if no nonzero element
of~$\su(m)$ commutes with both~$j_{Z_1}$ and~$j_{Z_2}$\,.
\endroster

\proclaim{3.2.5 Proposition}
Let $j,j':\h\cong\R^2\to\su(m)$ be two linear maps, and let~$g\subla\,,
g\sublap$ be the associated pair of Riemannian metrics on~$S^{2m+1}$,
or on~$B^{2m+2}$, as above.
If~$j\sim j'$ then the Riemannian manifolds $(S^{2m+1},g\subla)$
and $(S^{2m+1},g\sublap)$ are isospectral, and $(B^{2m+2},g\subla)$
and $(B^{2m+2},g\sublap)$ are Dirichlet and Neumann isospectral.
If~$j$ and~$j'$ are not equivalent and if at least one of them is
generic, then in both of these pairs the two manifolds are not isometric.
\endproclaim

\demo{Proof}
We will prove the nonisometry statement in Section~4.
To prove isospectrality we use Theorem~1.6.
Fix an arbitrary $\mu\in\LL^*$. We have to show that there exists
an isometry~$F_\mu$ of~$g_0$ which commutes with the action of~$T$
and satisfies $\mu\circ\la=F_\mu^*(\mu\circ\lap)$.
Let $Z\in\h$ be the vector corresponding to $\mu\in\LL^*\subset\h^*$
under the canonical identification of~$\h$ with~$\h^*$ associated
with the basis $\{Z_1\,,Z_2\}$.
Choose $A_Z\in\SU(m)\subset\SO(2m)$ as in Definition~3.2.4(i),
and let $F_\mu:=(A_Z\,,\Id)\in\SO(2m+2)$.
Then~$F_\mu$ is an isometry of~$g_0$ and satisfies
$$\align
\bigl(F_\mu^*(\mu\circ\lap)\bigr)_{(p,q)}(X,U)&=
(\mu\circ\lap)_{(A_Zp,q)}(A_ZX,U)\\
&=|A_Zp|^2\<j'_ZA_Zp,A_ZX\>-\<A_ZX,iA_Zp\>\<j'_ZA_Zp,iA_Zp\>\\
&=|p|^2\<A_Z\inv j'_ZA_Zp,X\>-\<X,ip\>\<A_Z\inv j'_ZA_Zp,ip\>\\
&=|p|^2\<j_Zp,X\>-\<X,ip\>\<j_Zp,ip\>=(\mu\circ\la)_{(p,q)}(X,U),
\endalign
$$
as desired.
\qed\enddemo

The following result shows that there are in fact many examples to
Proposition~3.2.5.

\proclaim{3.2.6 Proposition}
Let $m\ge3$ and $\{Z_1\,,Z_2\}$ be the standard basis of~$\R^2$.
\roster
\item"(i)"
There exists a nonempty Zariski open subset~$\Cal U$
of the space~$\J$ of linear maps $j:\R^2\to\su(m)$ such that for each
$j\in\Cal U$ there is a continuous family~$j(t)$ in~$\J$, defined on some
open interval around~$t=0$, such that $j(0)=j$ and:\newline
\qquad\qquad{\rm1.)} The maps~$j(t)$ are pairwise isospectral in the sense
of~{\rm3.2.4(i)}.\newline
\qquad\qquad{\rm2.)} The function $t\mapsto\|j_{Z_1}(t)^2+j_{Z_2}(t)^2\|^2=
  \tr\bigl((j_{Z_1}(t)^2+j_{Z_2}(t)^2)^2\bigr)$ is not constant in~$t$
  in any interval around zero. In particular, the maps~$j(t)$ are not
  pairwise equivalent in the sense of~{\rm3.2.4(ii)}.\newline
  \qquad\qquad{\rm3.)} The maps~$j(t)$ are generic in the sense
  of~{\rm3.2.4(iii)}.
\item"(ii)"
For $m=3$, an explicit example of an isospectral family~$j(t):\R^2\to
\su(3)$ with
$\|j_{Z_1}(t)^2+j_{Z_2}(t)^2\|^2\ne
\text{\rm const}$ is given by
$$j_{Z_1}(t):=\left(\smallmatrix -i&0&0\\0&0&0\\0&0&i\endsmallmatrix\right),
\qquad j_{Z_2}(t):=\left(\smallmatrix 0&\cos t&\sqrt2\sin t\\-\cos t
 &0&\cos t\\-\sqrt2\sin t&-\cos t&0\endsmallmatrix\right).
$$
The $j(t)$ are pairwise isospectral since
$\det\bigl(\la\Id-(sj_{Z_1}(t)+uj_{Z_2}(t))\bigr)=\la^3+(s^2+2u^2)\la$
is independent of~$t$. However, $\|j_{Z_1}(t)^2+j_{Z_2}(t)^2\|^2=
14+4\,{\sin^2}t$ is nonconstant in~$t$.
The map~$j(t)$ is generic in the sense of~{\rm3.2.4(iii)} if and only
if~$\cos t\ne0$.
\endroster\endproclaim

\demo{Proof}
(i)
The same statement, but without requirement~3.), was shown in
\cite{Sch3, Proposition~3.6}. Since the section of two nonempty Zariski
open sets is again such a set, it just remains to show that the
subset~$\Cal O\subset\Cal J$ of those elements which are generic in
the sense of~3.2.4(iii) is nonempty and Zariski open.
Note that the map $j\in\J$ given by
$$j_{Z_1}:=\left(\smallmatrix i\alpha_1&&&&\\&i\alpha_2&&&\\&&i\alpha_3&&\\
&&&\ddots&\\&&&&i\alpha_m\endsmallmatrix\right),
\qquad j_{Z_2}:=\left(\smallmatrix 0&-1&&&\\1&0&-1&&\\&1&\ddots&\ddots&\\
&&\ddots&0&-1\\&&&1&0\endsmallmatrix\right)
$$
with pairwise different
$\alpha_1\,,\,\ldots\,\alpha_m\in\R$ satisfying
$\alpha_1+\ldots+\alpha_m=0$ is an element of~$\Cal O$; thus~$\Cal O$
is nonempty.
To see that~$\Cal O$ is Zariski open, note that it is equal to the
set of those $j\in\J$ for which the map $F_j:\su(m)\ni\tau\mapsto
([j_{Z_1}\,,\tau],[j_{Z_2}\,,\tau])\in\su(m)\oplus\su(m)$
has maximal rank $r_m:=\dimm(\su(m))$.
But the latter condition can be expressed as the nonvanishing of a
certain polynomial in the coefficients of~$j$, namely, the sum of the
squared determinants of the $(r_m\times r_m)$@-minors of a matrix
representation of~$F_j$\,.

(ii) This can be checked by straightforward calculation.
\qed\enddemo

\proclaim{3.2.7 Corollary}
For every $m\ge3$ there exist continuous families of isospectral
metrics on~$S^{2m+1}$ and continuous families of Dirichlet
and Neumann isospectral metrics on~$B^{2m+2}$.
In particular, there exist such families on~$S^7$, resp\. on~$B^8$.
An explicit example is given by the metrics~$g_{\la(t)}$ associated
with the family~$j(t)$ from~{\rm3.2.6(ii)}.
\endproclaim

\remark{Remark}
Carolyn Gordon \cite{Go3} has previously given continuous families of
isospectral metrics on each $S^{n\ge8}$ and $B^{n\ge9}$ using a related
construction (see~3.4(i)).
Those were the first examples
of continuous isospectral families of metrics on balls and spheres.
\endremark

\subheading{3.3 Example: Isospectral pairs of metrics on
$\boldkey S^\boldkey5$ and on $\boldkey B^\boldkey6$}

\noindent
\subsubhead 3.3.1 Notation\endsubsubhead
In the context of Notation~3.1 we consider now the case $m=2$.

(i) We fix a realization of the Hopf projection $P:S^3\to S^2_{1/2}\subset
\R^3$ in coordinates, say
$$P:(\alpha,\beta,\gamma,\delta)\mapsto\bigl(\tfrac12(\alpha^2+\beta^2
  -\gamma^2-\delta^2),\,\alpha\gamma+\beta\delta,\,\alpha\delta-\beta\gamma
  \bigr).
$$
We extend~$P$ to a smooth map from~$\R^4\cong\C^2$ to~$\R^3$, defined
by the same formula (note that~$P$ will map~$S^3_a$ to~$S^2_{a^2/2}$
for each radius $a\ge0$).

(ii) Let $\Sym_0(\R^3)$ denote the space of symmetric traceless real
$(3\times3)$@-matrices.
For each linear map $c:\h\cong\R^2\to\Sym_0(\R^3)$ we define a $\h$@-valued
$1$@-form $\la=(\la^1,\la^2)$ on $\R^6\cong\C^2\oplus\C$ by letting
$$\la^k_{(p,q)}(X,U)=\<c_{Z_k}P(p)\times P(p),P_{*p}(X)\>\tag{5}
$$
for $k=1,2$ and all $(X,U)\in T_p\R^{2m}\oplus T_q\R^2$,
where $\scp$ denotes the standard euclidean inner product on~$\R^3$,
and~$\times$ denotes the vector product in~$\R^3$.
By restriction we obtain a smooth $\h$@-valued $1$@-form~$\la$ on the
unit sphere~$S^5$, respectively on the unit ball~$B^6$.

\subsubhead 3.3.2 Remark\endsubsubhead
We observe that~$\la$ is admissible
with respect to the action~$\rho$ of~$T$ on $S^5\subset\C^2\oplus\C$,
resp\. $B^6\subset\C^2\oplus\C$.
Invariance of~$\la$ under the action of~$T$ is immediate since
$P(e^{ia}p)=P(p)$ for all $a\in\R$, $p\in\C^2$. By the same reason we have
$P_{*p}(ip)=0$, which implies that~$\la$ vanishes on the
spaces~$\h_{(p,q)}$\,.

\subsubhead 3.3.3 Definition\endsubsubhead
Let~$g\subla$ be the metric associated with~$\la$ and the standard
metric~$g_0$ on~$S^5$, resp.~$B^6$.
Thus for each linear map $c:\R^2\to\Sym_0(\R^3)$ we have an associated
Riemannian metric~$g\subla$ on~$S^5$, respectively on $B^6$,
via the corresponding $1$@-form~$\la$ as defined in~\thetag{5}.

\subsubhead 3.3.4 Definition\endsubsubhead
Let $j,j':\h\cong\R^2\to \Sym_0(\R^3)$ be two linear maps.
\roster
\item"(i)"
We call $c$ and~$c'$ {\it isospectral}, denoted $c\sim c'$,
if for each $Z\in\h$ there exists $E_Z\in\SO(3)$ such that
$c'_Z=E_Zc_ZE_Z\inv$\,.
\item"(ii)"
We call $c$ and~$c'$ {\it equivalent}, denoted $c\cong c'$,
if there exists $E\in\O(3)$
and $\Psi\in\Cal E$ (see Notation~3.1(ii))
such that $c'_Z=Ec_{\Psi(Z)}E\inv$ for all $Z\in\h$.
\item"(iii)" We say $c$ is {\it generic\/} if no nonzero element
of~$\so(3)$ commutes with both~$c_{Z_1}$ and~$c_{Z_2}$\,.
\endroster

\proclaim{3.3.5 Proposition}
Let $c,c':\h\cong\R^2\to \Sym_0(\R^3)$ be two linear maps, and
let~$g\subla\,, g\sublap$ be the associated pair of Riemannian metrics
on~$S^5$, or on~$B^6$, as above.
If $c\sim c'$ then the Riemannian manifolds $(S^5,g\subla)$
and $(S^5,g\sublap)$ are isospectral, and $(B^6,g\subla)$
and $(B^6,g\sublap)$ are Dirichlet and Neumann isospectral.
If $c$ and $c'$ are not equivalent and if at least one of them is generic,
then in both of these pairs the two manifolds are not isometric.
\endproclaim

\demo{Proof}
Again, we postpone the proof of the nonisometry statement to Section~4,
and we use Theorem~1.6 to prove isospectrality.
Fix $\mu\in\LL^*$
and let~$Z$ be the dual vector in~$\h$ as in the proof of~3.2.5.
Choose $E_Z\in\SO(3)$ as in Definition~3.3.4(i),
and choose $A_Z\in\SU(2)\subset\SO(4)$ such that for the Hopf
projection $P:S^3\to S^2_{1/2}$ from 3.3.1(i) we have $P\circ A_Z
=E_Z\circ P$. Let $F_\mu:=(A_Z\,,\Id)\in\SO(6)$.
Then~$F_\mu$ is an isometry of~$g_0$ and satisfies
$$\align
\bigl(F_\mu^*(\mu\circ\lap)\bigr)_{(p,q)}&(X,U)=
(\mu\circ\lap)_{(A_Zp,q)}(A_ZX,U)\\
&=\<c'_ZP(A_Zp)\times P(A_Zp),P_{*A_Zp}(A_ZX)\>\\
&=\<c'_ZE_ZP(p)\times E_ZP(p),E_ZP_{*p}(X)\>\\
&=\<E_Z\inv c'_ZE_ZP(p)\times P(p),P_{*p}(X)\>\\
&=\<c_ZP(p)\times P(p),P_{*p}(X)\>
=(\mu\circ\la)_{(p,q)}(X,U),
\endalign
$$
as desired.
\qed\enddemo

\proclaim{3.3.6 Proposition}
There exist pairs of linear maps $c,c':\R^2\to\Sym_0(\R^3)$ such
that~$c$ and~$c'$ are isospectral in the sense of {\rm3.3.4(i)},
not equivalent in the sense of {\rm3.3.4(ii)}, and both generic in the
sense of {\rm3.3.4(iii)}.
An example of such a pair is given by
$$
\hskip29pt c_{Z_1}=c'_{Z_1}=\left(\smallmatrix -1&0&0\\0&0&0\\0&0&1
\endsmallmatrix\right),\ \
c_{Z_2}=\left(\smallmatrix 0&1&0\\1&0&1\\0&1&0\endsmallmatrix\right),\ \
c'_{Z_2}=\left(\smallmatrix 0&0&\sqrt2\\0&0&0\\ \sqrt2&0&0\endsmallmatrix
\right),
$$
where $\{Z_1\,,Z_2\}$ is the standard basis of~$\R^2$.
\endproclaim

\demo{Proof}
One easily checks that for each fixed pair $s,u\in\R$, the
characteristic polynomials of $sc_{Z_1}+uc_{Z_2}$ and $sc'_{Z_1}
+uc'_{Z_2}$ are equal (namely, to $\la^3+(s^2+2u^2)\la$),
which implies isospectrality.
That~$c$ and $c'$ are not equivalent can be seen from the fact
that $\|c_{Z_1}^{\,2}+c_{Z_2}^{\,2}\|^2=14$,
while $\|c^{\prime\,2}_{Z_1}
+c^{\prime\,2}_{Z_2}\|^2=18$. The maps are generic in the sense
of 3.3.4(iii) because not even the map $c_{Z_1}=c'_{Z_1}$
alone commutes with any nonzero element of~$\so(3)$.
\qed\enddemo

\proclaim{3.3.7 Corollary}
There exist nontrivial pairs of isospectral metrics on~$S^5$,
and there exist nontrivial pairs of Dirichlet and Neumann isospectral
metrics on~$B^6$.
\endproclaim

\remark{Remark}
Our above pairs of isospectral $5$@-spheres constitute the lowest
dimensional examples of isospectral spheres which have been constructed
so far. The analogous statement holds for our isospectral
$6$@-dimensional balls.
\endremark

\subheading{3.4 Survey of related examples}
As a complement to our above new examples of isospectral spheres
and balls we now give a short survey of related examples.
Among them, (ii) and (v) are new; (i), (iii), and (iv) are
already known, but were constructed first in slightly different
settings.

Although we will not present any proofs here, we
note that in each of the examples below
the Riemannian manifolds can be described within the
setting of~1.5 and proven to be isospectral by Theorem~1.6.
Moreover, for each of them there is a nonisometry proof using the
general results from Section~2 and following the
lines of Section~4 below, where we will prove nonisometry
for the above examples~3.2\,/\,3.3.
Since the formulas defining~$\la,\lap$ in the examples below
(except~(v)) are less
complicated than those in~3.2 or~3.3, the corresponding
computations are simpler than
those in Section~4. Also note that in some of the examples the situation
is simplified by the fact that $\hat M=M$ (in~(iii) and~(v)),
or that~$\Om_0=0$ (in (i)--(iv)). For example, in~(iii) below,
the curvature form~$\Om\subla$ equals~$d\bar\la$ on $M/T=S^{2m-1}$,
where $d\bar\la^k$ is of the form $(X,Y)\mapsto 2\<j_{Z_k}X,Y\>$.

The notation we use is similar to the one we used above.
In particular, $T$~again denotes the torus $\R^2/(2\pi\Z\times
2\pi\Z)$, and $\{Z_1\,,Z_2\}$ is the standard basis of its Lie
algebra $\h\cong\R^2$.
On all manifolds~$M$ which we consider below there is a canonical
standard metric which will in each case play the role of~$g_0$\,.
When we call two manifolds with boundary isospectral we always mean
Dirichlet and Neumann isospectral.

\smallskip

{\bf (i)} {\it Continuous families of isospectral metrics on~$S^{n\ge8}$
and on~$B^{n\ge9}$ \cite{Go3}.}
\newline
Let $M=S^{m+3}\subset\R^m\oplus\C\oplus\C$, resp\. $M=B^{m+4}\subset
\R^m\oplus\C\oplus\C$. Let the action of~$T$ on~$M$ be induced by its
canonical action on the $\C\oplus\C$ component (generated by multiplication
with~$i$ on each summand). For each linear map $j:\R^2\to\so(m)$
we define a $\h\cong\R^2$@-valued $1$@-form on $\R^m\oplus\C\oplus\C
\cong\R^m\oplus\R^4$ (and hence on~$M$) by letting
$$\la^k_{(p,q)}(X,U)=\<j_{Z_k}p,X\>
$$
for $k=1,2$ and all $(X,U)\in T_p\R^m\oplus T_q\R^4$.
The following conditions on a pair $j,j'$ imply that the associated
Riemannian manifolds $(M,g\subla)$, $(M,g\sublap)$ are isospectral
and not isometric:
\roster
\item"1.)" $j$ and~$j'$ are isospectral: For each $Z\in\h$ there exists
$A_Z\in\O(m)$ such that $j'_Z=A_Zj_ZA_Z\inv$.
\item"2.)" $j$ and $j'$ are nonequivalent: There is no $A\in\O(m)$
and no $\Psi\in\tilde\Cal E$ such that $j'_Z=Aj_{\Psi(Z)}
A\inv$ for all $Z\in\h$;
here~$\tilde\Cal E$ denotes the set of the eight automorphisms of~$\R^2$
which preserve the set $\{\pm Z_1\,,\pm Z_2\}$.
\item"3.)" At least one of~$j,j'$ is generic; $j$~is called generic if
no nonzero element of~$\so(m)$ commutes with both~$j_{Z_1}$ and~$j_{Z_2}$\,.
\endroster
It is known~\cite{GW2} that for $m\ge5$ there are continuous families
-- even multipara\-meter families --~$j(t)$ whose elements pairwise satisfy
these conditions. The associated manifolds $(M,g_{\la(t)})$ are
Carolyn Gordon's examples of isospectral deformations on spheres and balls.

\smallskip

{\bf (ii)} {\it Isospectral pairs of metrics on~$S^6$ and on~$B^7$.}
\newline
This is just a slight modification of~Example~3.3 which does not involve
the Hopf projection~$P$ anymore; in turn, an additional dimension is needed
for the construction.
Let $M=S^6\subset\R^3\oplus\C\oplus\C$, resp\. $M=B^7\subset
\R^3\oplus\C\oplus\C$. Let~$T$ act canonically on the $\C\oplus\C$ component
as in~(i). For each linear map $c:\R^2\to\Sym_0(\R^3)$
we define a $\h\cong\R^2$@-valued $1$@-form on $\R^3\oplus\C\oplus\C
\cong\R^3\oplus\R^4$ (and hence on~$M$) by letting
$$\la^k_{(p,q)}(X,U)=\<c_{Z_k}p\times p,X\>
$$
for $k=1,2$ and all $(X,U)\in T_p\R^3\oplus T_q\R^4$.
The conditions on a pair~$c,c'$ which cause $(M,g\subla)$ and $(M,g\sublap)$
to be isospectral and not isometric are the same as in~3.3, except that
in the nonequivalence condition the group~$\Cal E$ has to be replaced
by the larger group~$\tilde\Cal E$ as in~(i).

The specific pair~$c,c'$ given in~3.3.6 satisfies these
conditions and thus yields a pair of isospectral metrics~$g\subla\,,
g\sublap$ on~$S^6$ (resp\. on~$B^7$). Moreover, one can compute the
loci~$N,N'$ of the maximal scalar curvature in~$(S^6,g\subla)$
and~$(S^6,g\sublap)$. It turns out that~$N'$ contains a $4$@-sphere,
while~$N$ is a union of $2$@- and $3$@-spheres.

\smallskip

{\bf (iii)} {\it Continuous families of isospectral metrics on~$S^{m-1\ge4}
\times T^2$ \cite{GGSWW} and on $B^{m\ge5}\times T^2$ \cite{GW2};
pairs of isospectral metrics on $S^2\times T^2$ \cite{Sch3} and on
$B^3\times T^2$.}
\newline
For $a,b,c>0$ with $a^2+b^2+c^2\le1$ we define
$M_{a,b,c}:=\{(p,u,v)\in\R^m\oplus\C\oplus\C\mid|p|^2=a^2,|u|^2=b^2,
|v|^2=c^2\}$. Then~$M_{a,b,c}$ is diffeomorphic to $S^{m-1}\times T^2$
and is a submanifold of $B^{m+4}\subset\R^m\oplus\C\oplus\C$.
The restrictions of two metrics $g\subla\,,g\sublap$ from~(i)
to one of these submanifolds are again isospectral and nonisometric
under the same conditions as in~(i) on the underlying pair of maps~$j,j'$.
In particular, one obtains continuous isospectral families~$g_{\la(t)}$
on $S^{m-1}\times T^2$ for $m\ge5$. These are just the examples constructed
in \cite{GGSWW}.

Similarly, for $m=3$ the metrics $g\subla\,,g\sublap$ from~(ii) restrict
to isospectral nonisometric metrics on $M_{a,b,c}$ which is now diffeomorphic
to $S^2\times T^2$; these pairs were first constructed in~\cite{Sch3}.

In both cases we can replace the condition $|p|^2=a^2$ by $|p|^2\le a^2$
in the definition of~$M_{a,b,c}$ which then becomes diffeomorphic
to $B^m\times T^2$. We obtain continuous families of isospectral
metrics~$g_{\la(t)}$
on $B^m\times T^2$ for $m\ge5$ \cite{GW2} and pairs of such metrics
on $B^3\times T^2$.

\smallskip

{\bf (iv)} {\it Isospectral pairs of metrics on $S^2\times S^3$
\cite{Ba}.}
\newline
Instead of considering the submanifolds $M_{a,b,c}\approx S^2\times T^2$
of~$B^7\subset\R^3\oplus\C\oplus\C$ as in the middle part of~(iii) above,
we now consider $M_{a,b}:=\{(p,q)\in\R^3\oplus\R^4\mid |p|^2=a^2,
|q|^2=b^2\}$ with $a,b>0$. Then $M_{a,b}$ is diffeomorphic to $S^2\times
S^3$, and again the metrics $g\subla\,,g\sublap$ from~(ii) restrict
to isospectral nonisometric metrics on~$M_{a,b}$\,.

\smallskip

{\bf (v)} {\it Continuous isospectral families on $S^{2m-1\ge5}\times S^1$
and on $B^{2m\ge6}\times S^1$; pairs of isospectral metrics on
$S^3\times S^1$ and on $B^4\times S^1$.}
\newline
We return to the context of our examples~3.2\,/\,3.3 above and consider
the submanifolds $M_{a,b}:=\{(p,q)\in\C^m\oplus\C\mid|p|^2=a^2,
|q|^2=b^2\}$ of~$B^{2m}$. For $a,b>0$ the manifold~$M_{a,b}$ is
diffeomorphic to $S^{2m-1}\times S^1$. If $j,j':\R^2\to\su(m)$
satisfy the assumptions of Proposition~3.2.5 then the associated pair
of metrics $g\subla\,,g\sublap$ on~$B^{2m}$ restricts to a pair of
isospectral nonisometric metrics on~$M_{a,b}$ as well. In particular,
for $m\ge3$ we obtain continuous families~$g_{\la(t)}$ of such metrics
as in~3.2.6\,/\,3.2.7.

Similarly, for $m=2$ each pair $c,c':\R^2\to\Sym_0(\R^3)$ satisfying
the assumptions of Proposition~3.3.5 also yields a pair of isospectral
nonisometric metrics on the manifold $S^{2m-1}\times S^1=S^3\times S^1$.

In both cases we can again replace the condition $|p|^2=a^2$
by $|p|^2\le a^2$ in the definition of~$M_{a,b}$ which then
becomes diffeomorphic to $B^{2m}\times S^1$. We obtain continuous
families of isospectral metrics~$g_{\la(t)}$ on $B^{2m}\times S^1$
for $m\ge3$, and pairs of such metrics on $B^4\times S^1$.

\bigskip

\heading\S4 Nonisometry of the examples from 3.2\,/\,3.3
\endheading

Our strategy for proving the nonisometry statements of Proposition
3.2.5\,/\,3.3.5 consists in showing that the nonequivalence condition
on~$j,j'$ (resp\. on~$c,c'$) implies condition~\thetag{N} of Proposition~2.4,
while the genericity condition on the pair of maps implies
Property~\thetag{G} for the associated curvature forms. See Proposition~4.3
below for these assertions. The desired nonisometry statements then
follow immediately from Proposition~2.4.

\subheading{4.1 Notation and Remarks}
Throughout this section
we denote by~$(M,g_0)$ either $S^{2m+1}\subset\R^{2m}
\oplus\R^2$ or $B^{2m+2}\subset\R^{2m}\oplus\R^2$, endowed with the
standard metric. We consider the
action~$\rho$ of $T=\R^2/(2\pi\Z\times2\pi\Z)$ on~$M$ which was defined
in~3.1(iii).
\roster 
\item"(i)" 
For $a,b\ge0$ let $M_{a,b}:=\{(p,q)\in\R^{2m}\oplus\R^2\mid |p|^2=a^2, 
|q|^2=b^2\}$. Thus~$S^{2m+1}$, resp\. $B^{2m+2}$, is the disjoint union of 
the submanifolds~$M_{a,b}$ with $a^2+b^2=1$, resp\. $a^2+b^2\le1$. 
Note that in either case, $\hat M$~is the union of those~$M_{a,b}\subset M$
with $a,b>0$. Each~$M_{a,b}$ is obviously $T$@-invariant.
\item"(ii)" 
For $M_{a,b}\subset\hat M$ the manifold $(M_{a,b}/T,g_0^T)$ is isometric to 
$(\C P^{m-1},a^2\gfs)$, where~$\gfs$ denotes the Fubini-Study metric.
\item"(iii)"
The first component of the $\h\cong\R^2$@-valued form which
the curvature form~$\Om_0$ induces
on $(M_{a,b}/T,g_0^T)$ $=$ $(\C P^{m-1},a^2\gfs)$ is a scalar multiple of
the standard K\"ahler form on $\C P^{m-1}$, and the second component is zero.
In fact,
for $(X,Z),(Y,W)\in T_{(p,q)}M_{a,b}=T_pS_a^{2m-1}\oplus T_qS_b^1$
we have $\om_0^1(X,Z)=\<X,ip\>/a^2$, $\om_0^2(X,Z)=\<Z,iq\>/b^2$, hence
$d\om_0^1((X,Z),(Y,W))=2\<iX,Y\>/a^2$ and
$d\om_0^2((X,Z),(Y,W))$ $=$ $2\<iZ,W\>/b^2=0$,
where the last equation follows from the fact that $T_qS_b^1$ is
one-dimensional. The statement now follows because~$\Om_0$ is induced
by~$d\om_0$\,.
\item"(iv)"
Any isometry of $(\C P^{m-1},a^2\gfs)=S^{2m-1}_a/S^1$ is induced by
some $\C$@-linear or $\C$@-antilinear isometry $A\in\SU(m)\cup
\SU(m)\circ Q$ (see 3.2.4) of~$\C^m$.
\endroster

\proclaim{4.2 Lemma}
If $F\in\Autm$ then $F$~preserves each of the submanifolds~$M_{a,b}$
of~$\hat M$. Moreover, $\Psi_F\in\Cal E$ {\rm(}see Notation~{\rm3.1(ii))}.
In particular, for the group~$\Cal D$ defined in~{\rm2.1(ii)} we have
$\Cal D\subset\Cal E$.
\endproclaim

\demo{Proof}
Note that the $T$@-orbit through $(p,q)\in M_{a,b}\subset\hat M$, 
endowed with the metric induced by~$g_0$\,, 
is a rectangular torus with side lengths $2\pi a$ and $2\pi b$. 
It follows that~$F$ preserves the sets $M_{a,b}\cup M_{b,a}\subset\hat M$. 
We have to show that~$F$ cannot switch the components $M_{a,b}$ 
and~$M_{b,a}$ if $a\ne b$. If it did then it would induce an isometry 
from $(M_{a,b}/T,g_0^T)$ to $(M_{b,a}/T,g_0^T)$. 
But by~4.1(ii) these manifolds have different volume.  
Thus~$F$ preserves~$M_{a,b}$\,.
Now choose $0<a<b$ such that $M_{a,b}\subset\hat M$,
and let $(p,q)\in M_{a,b}$\,. Since~$F$ preserves
$M_{a,b}$\,, both of the orbits $T\cdot(p,q)$ and $T\cdot F(p,q)$ 
are rectangular tori on which the shortest closed loops have length 
$2\pi a$ and are precisely the flow lines of~$Z_1^*$\,. This, together 
with $Z_1^*\perp Z_2^*$ and $\|Z^*_2\restr{(p,q)}\|=\|Z^*_2\restr{F(p,q)}\| 
=b$, implies that $F_*(Z_1^*)\in\{\pm Z_1^*\}$ and $F_*(Z_2^*)\in 
\{\pm Z_2^*\}$; hence $\Psi_F\in\Cal E$.
Now $\Cal D\subset\Cal E$ follows from the definition of~$\Cal D$.
\qed\enddemo 

In view of Proposition~2.4, the following result implies immediately
the nonisometry statements of Proposition 3.2.5\,/\,3.3.5\,:

\proclaim{4.3 Proposition}
Let $\la,\lap$ be of the type defined in~\thetag{4} {\rm(}resp\.
in~\thetag{5}{\rm)},
and let $j,j':\R^2\to\su(m)$ {\rm(}resp\. $c,c':\R^2\to\Sym_0(\R^3)${\rm)}
be the pair of linear maps with which $\la,\lap$ are associated.
\roster
\item"(i)"
If~$j$ and~$j'$ {\rm(}resp.~$c$ and~$c'${\rm)} are not equivalent in the
sense of {\rm3.2.4(ii) (}resp.~{\rm3.3.4(ii))}, then~$\Om\subla$
and~$\Om\sublap$ satisfy condition~\thetag{N} of Proposition~{\rm2.4}.
\item"(ii)"
If $j$ {\rm(}resp\. $c${\rm)} is generic in the sense of
{\rm3.2.4(iii) (}resp\. {\rm 3.3.4(iii))}, then~$\Om\subla$
has Property~\thetag{G}.
\endroster
\endproclaim

\demo{Proof}
We choose one of the submanifolds $M_{a,b}=:L$ of~$\hat M$
and denote by $\om\subla^L\,,\Om\subla^L$ the forms induced by~$\om\subla\,,
\Om\subla$ on $L\subset\hat M$ and $L/T\subset\hat M/T$, respectively.

(i)
Suppose that condition~\thetag{N} were not satisfied.
Let $\Psi\in\Cal D\subset\Cal E$ (recall~4.2) and $\bar F\in\Autbm$
such that $\Om\subla=\Psi\circ\bar F^*\Om\sublap$\,.
Lemma~4.2 implies that~$\bar F$ preserves $L/T$; hence
$\Om\subla^L=\Psi\circ\bar F^*\Om\sublap^L$.
Recall from~2.1(vi) that $\Om\subla=\Om_0+d\bar\la$, $\Om\sublap=\Om_0+d
\bar\lap$.
Each closed $2$@-form on~$\C P^{m-1}$ is uniquely decomposable into an
exact component and a multiple of the K\"ahler form. Therefore, the
description of~$\Om_0^L$ given in~4.1(iii) now implies that
$$d\bar{\la}^L=\Psi\circ\bar F^*d\bar{\lap}{}^L\tag{6}
$$
and hence
$$d\la^L=\Psi\circ A^*d{\lap}{}^L\tag{7}
$$
for the map $A\in\SU(m)\cup\SU(m)\circ Q$
which induces the isometry~$\bar F$ of $(L/T,g_0^T)=(\C P^{m-1},
a^2\gfs)$ (recall~4.1(iv)).
Here $\la^L,\lap{}^L$ and $\bar\la^L,\bar{\lap}{}^L$
denote the forms which $\la,\lap$ and $\bar\la,\bar{\lap}$
induce on~$L$ and~$L/T$, respectively.
We are going to show that
$$j'_{\Psi(Z)}=Aj_ZA\inv\tag{8}
$$
for all $Z\in\h$ in the case of forms of type~\thetag{4}, and
$$c'_{\Phi(Z)}=Ec_ZE\inv\tag{9}
$$
for all $Z\in\h$ in the case of forms of type~\thetag{5}, where
$E:=\bar F\restr{L/T}\in\O(3)$ and $\Phi:=\det(E)\Psi\in\{\pm\Psi\}\subset
\Cal E$. This will contradict the assumed nonequivalence of~$j$ and~$j'$
(resp\. of~$c$ and~$c'$).

In the first case, we have
for all $(X,Z),(Y,W)\in T_{(p,q)}L=T_pS_a^{2m-1}\times
T_qS_b^1$ and for $k=1,2$:
$$\align\la^k(X,Z)&=a^2\<j_kp,X\>-\<j_kp,ip\>\<X,ip\>,\quad\text{hence}\\
d\la^k((X,Z),(Y,W))&=2a^2\<j_kX,Y\>-2\<j_kX,ip\>\<Y,ip\>+2\<j_kY,ip\>
 \<X,ip\>\\
&\quad-2\<j_kp,ip\>\<iX,Y\>\\
&=2a^2\<j_kX^h,Y^h\>-2\<j_kp,ip\>\<iX,Y\>,\endalign
$$
where $j_k:=j_{Z_k}$ and $X^h$~denotes the $g_0$@-orthogonal projection
of~$X$ to~$(ip)^\perp$. Let $\eps_k\in\{\pm1\}$ be such that
$\Psi(Z_k)=\eps_kZ_k$\,. Then one derives
from equation~\thetag{7}:
$$\split \eps_k(a^2\<j_k'(AX)^h,(AY)^h\>-\<j_k'Ap,iAp\>\<iAX,AY\>)\\
         =a^2\<j_kX^ h,Y^h\>-\<j_kp,ip\>\<iX,Y\>\endsplit
$$
for $k=1,2$ and all $X,Y\in T_pS_a^{2m-1}$.
Since $(AX)^h=AX^h$ and since~$A$ either commutes or anticommutes with~$i$,
we get, letting $\tau_k:=\eps_kA\inv j'_kA-j_k\in\su(m)$:
$$a^2\<\tau_kX^h,Y^h\>-\<\tau_kp,ip\>\<iX,Y\>=0.
$$
In particular, we obtain for all $p\in S_a^{2m-1}$ and all nonvanishing
$X\in\spann\{p,ip\}^\perp$, by letting $Y:=iX$:
$$\<\tau_kX,iX\>/|X|^2=\<\tau_kp,ip\>/|p|^2.
$$
The hermitian map~$i\tau_k$ must therefore be a scalar multiple of the
identity,
and~$\tau_k$ be a scalar multiple of $i\Id$. This implies $\tau_k=0$
because of
$\tau_k\in\su(m)$;
equation~\thetag{8} now follows.

In the second case (namely, with $\la,\lap$ of type~\thetag{5})
we have
$$\bar\la^k(X)=a^3\<c_kx\times x,X\>
$$
for $k=1,2$ and all $X\in T_x(L/T)=T_xS^2_{a/2}$\,, where $c_k:=c_{Z_k}$\,.
The factor~$a^3$ is due to the fact that if $\pi:S_a^3\to S^2_{a/2}$
denotes the Riemannian submersion, then for the projection~$P$ as defined
in~3.3.1 we have $|P(p)|=a|\pi(p)|$. We obtain
$$\align d\bar\la^k(X,Y)/a^3&=\<c_kX\times x+c_kx\times X,Y\>
  -\<c_kY\times x+c_kx\times Y,X\>\\
  &=\<c_kX\times x,Y\>+2\<c_kx\times X,Y\>-\<c_kY\times x,X\>\\
  &=\<c_kX\times x,Y\>+2\<c_kx\times X,Y\>+\<Y\times c_kx,X\>+\<Y\times x,
  c_kX\>\\
  &=3\<c_kx\times X,Y\>.\endalign
$$
Note that in the third equation we have used $\tr(c_k)=0$.
Equation~\thetag{6} now implies for $\bar F\restr{S^2_{a/2}}=E\in\O(3)$:
$$\eps_k\<c_k'Ex\times EX,EY\>=\<c_kx\times X,Y\>
$$
for $k=1,2$ and all $X,Y\in T_xS^2_{a/2}$\,, where~$\eps_k$ is defined
as above. Letting $\tau_k:=\eps_k\det(E)E\inv c'_kE-c_k\in\Sym_0(\R^3)$
we obtain
$$\<\tau_kx\times X,Y\>=0
$$
for all $X,Y\perp x$, which implies $\tau_kx\perp x$ for all
$x\in S^2_{a/2}$\,. Since~$\tau_k$ is symmetric, it must be zero;
equation~\thetag{9} now follows.

(ii)
Suppose to the contrary that~$\la$ does not have Property~\thetag{G};
i.e., there is a $1$@-parameter family $\bar F_t\in\Autbm$ such
that $\bar F_t^*\Om\subla\equiv\Om\subla$\,.
Proceeding as in the proof of~(i)
(in the special situation $\Psi=\Id$ and $\la=\lap$)
we now obtain $1$@-parameter families $A_t\in\SU(m)$ (resp\. $E_t\in\SO(3)$)
preserving~$d\la^L$ (resp.~$d\bar\la^L$), and we derive
that $j_Z\equiv A_tj_ZA_t\inv$ (resp\. $c_Z\equiv E_tc_ZE_t\inv$)
for each $Z\in\h$ (cp.~\thetag{8},~\thetag{9}).
But this contradicts the genericity assumption made on~$j$ (resp.~$c$).
\qed\enddemo

\bigskip

\heading\S5 Isospectral metrics on spheres and balls which are equal
to the standard metric on large subsets\endheading

\noindent
The main idea of this section is to simultaneously multiply~$\la,\lap$
by some smooth scalar function~$f$ on the sphere (resp\. the ball)
which has small support but is chosen such that both isospectrality
and nonisometry of the associated metrics $g_{f\la}$ and $g_{f\lap}$
continue to hold.

The following observation is trivial.

\subheading{5.1 Remark}
In the context of Theorem~1.6, let $\la,\lap$ be two admissible $\h$@-valued
$1$@-forms on~$M$ which satisfy~\thetag{$*$}.
Let $f\in C^\infty(M)$ be any function which is invariant under~$T$ and
under each of the isometries~$F_\mu$ occurring in~\thetag{$*$}.
Then $f\la$ and $f\lap$ are again admissible and satisfy~\thetag{$*$}.
In particular, $(M,g_{f\la})$ and $(M,g_{f\lap})$ are isospectral.

\smallskip

\proclaim{5.2 Proposition}
Let $(M,g_0)=S^{2m+1}\subset\R^{2m}\oplus\R^2$, resp\. $(M,g_0)
=B^{2m+2}\subset \R^{2m}\oplus\R^2$. Choose any $\phi\in C^\infty([0,1]^2)$
which does not vanish identically; in case $M=S^{2m+1}$ we assume that~$\phi$
does not vanish identically on the set $\{(s,1-s)\mid s\in[0,1]\}$.
Let $f\in C^\infty(M)$ be defined by
$f(p,q):=\phi(|p|^2,|q|^2)$ for $(p,q)\in M\subset\R^{2m}\oplus\R^2$.
Then under the assumptions of Proposition~{\rm3.2.5 (resp.~3.3.5)}
on the pair~$j,j'$ {\rm(}resp.~$c,c'${\rm)}
which defines~$\la,\lap$ as in~\thetag{4} {\rm(}resp.~\thetag{5}{\rm)},
the manifolds
$(M,g_{f\la})$ and $(M,g_{f\lap})$ are isospectral and not isometric.
\endproclaim

\demo{Proof}
The function~$f$ is obviously invariant under the action of~$T=\R^2/\LL$
as defined in~3.1(iii). Moreover, $f$~is invariant under those~$F_\mu$
occurring in the isospectrality proofs of Proposition 3.2.5\,/\,3.3.5:
Recall that the~$F_\mu$ used there were of the form $(A_Z\,,\Id)$ with
$A_Z\in\SO(2m)$. The isospectrality statement thus follows from
Remark~5.1 above.

For the nonisometry proof we must only slightly modify Proposition~4.3
and its proof.
In both statements of that proposition we replace each occurrence
of~$\Om\subla$ by~$\Om_{f\la}$\,, and similarly~$\Om\sublap$
by~$\Om_{f\lap}$\,. In the proof, we now choose
$L:=M_{a,b}\subset\hat M$ not arbitrarily,
but such that $C:=\phi(a^2,b^2)\ne0$.
This is possible by the assumptions made on~$\phi$.
The forms induced by $f\la$,~$f\lap$ on~$L$ are then equal to the
ones induced by~$C\la$ and~$C\lap$, respectively.
But replacing $\la^L,\la^{\prime\,L}$ by $C\la^L,C\la^{\prime\,L}$
does not afflict any of the subsequent arguments of the proof.
Thus both statements of Proposition~4.3, modified as described, remain true.
\qed\enddemo

\proclaim{5.3 Theorem}
Let $\eps>0$ be given, and let $m\ge2$.
Then there exist nonisometric pairs of isospectral metrics
on $S^{2m+1}$ whose volume element is the standard one
and which are equal to the round standard
metric outside a subset of volume smaller than~$\eps$.
For $m\ge3$ there are even continuous families of such metrics.
The~analogous statements hold for $B^{2m+2}$ {\rm(}with metrics
which are both Dirichlet and Neumann isospectral\,{\rm)}.
The mentioned sub\-set of small volume can be chosen as a tubular
neighborhood around any of the submanifolds $S_a^{2m-1}\times S_b^1$
with $0\le a,b\le1$ and
$a^2+b^2=1$ {\rm(}in the case of~$S^{2m+1}${\rm)}, resp\. $a^2+b^2\le1$
{\rm(}in the case of~$B^{2m+2}${\rm)}.
\endproclaim

\demo{Proof}
We can choose the function~$\phi\in C^\infty([0,1]^2)$
in Proposition~5.2 such that its support is contained in a sufficiently
small rectangle around~$(a^2,b^2)$.
Concerning the coincidence of the volume elements of~$g_0$ and~$g_{f\la}$
recall~1.5(iii).
\qed\enddemo

\enddocument